\documentclass[10pt,a4paper,reqno]{amsart}

\usepackage[centertags]{amsmath}

\usepackage{latexsym,amssymb}
\usepackage[ansinew]{inputenc}
\usepackage[T1]{fontenc}
\usepackage{times,mathptmx}
\usepackage{color}
\usepackage{graphicx,overpic}
\usepackage{hyperref}

\usepackage{amsfonts, times, mathptmx}
\usepackage{amssymb}
\usepackage{amsthm}
\usepackage{newlfont}
\usepackage[english]{babel}

\newcommand {\real} {\mathbb{R}}

\theoremstyle{plain}
\newtheorem{thm}{Theorem} [section]

\newtheorem{prop}[thm]{Proposition}

\theoremstyle{definition}

\theoremstyle{remark}

\begin {document}

\title[On the existence threshold for $p$-laplacian equations]{On the existence threshold for positive solutions of $p$-laplacian equations with a concave-convex nonlinearity}

\author {Fernando Charro, Enea Parini}

\date{\today}

\address{Department of Mathematics, The University of Texas at Austin, 1 University Station C1200, Austin, TX 78712} \email{fcharro@math.utexas.edu}

\address{Aix Marseille Universit\'{e}, CNRS, Centrale Marseille, I2M, UMR 7373, 13453 Marseille, France}

\email{enea.parini@univ-amu.fr}



\thanks{
Part of this work was carried over during a visit of the second author to the University of Texas at Austin. Both authors would like to thank this institution for the warm hospitality. The first author was partially supported by a MEC-Fulbright fellowship.
}

\begin {abstract} We study the following boundary value problem with a \emph{concave-convex} nonlinearity:
\begin{equation*}
\left\{
\begin{array}{r c l l}
-\Delta_p u  & = & \Lambda\,u^{q-1}+ u^{r-1} &  \textrm{in
}\Omega, \\  u & = & 0 & \textrm{on }\partial\Omega.
\end{array}\right.
\end{equation*} 
Here $\Omega \subset \real^n$ is a bounded domain and $1<q<p<r<p^*$. It is well known that there exists a number $\Lambda_{q,r}>0$ such that the problem admits at least two positive solutions for $0<\Lambda<\Lambda_{q,r}$, at least one positive solution for $\Lambda=\Lambda_{q,r}$, and no positive solution for $\Lambda > \Lambda_{q,r}$. We show that
\[ \lim_{q \to p} \Lambda_{q,r} = \lambda_1(p), \]
where $\lambda_1(p)$ is the first eigenvalue of the p-laplacian. It is worth noticing that $\lambda_1(p)$ is the threshold for existence/nonexistence of positive solutions to the above problem in the limit case $q=p$.

\end{abstract}

\maketitle

\section{Introduction}
Let $\Omega \subset \real^n$ be a bounded domain with Lipschitz boundary, and let us consider the boundary value problem
\begin{equation} \label{concaveconvex}
\left\{
\begin{array}{r c l l}
-\Delta_p u  & = & \Lambda\,u^{q-1}+ u^{r-1} &  \textrm{in
}\Omega, \\  u & = & 0 & \textrm{on }\partial\Omega,
\end{array}\right.
\end{equation}
where $1<q<p<r<p^*$ and $\Lambda > 0$. Due to the nature of the right-hand side, and in analogy with the case $p=2$, we will call Problem \eqref{concaveconvex} \emph{concave-convex}. It is well-known that there exists $\Lambda_{q,r}>0$ such that the problem admits at least two
positive solutions for every $\Lambda \in (0,\Lambda_{q,r})$, at least one positive solution for $\Lambda = \Lambda_{q,r}$, and no positive solution for $\Lambda > \Lambda_{q,r}$ (see \cite[Theorem 1.4]{garciaperalmanfredi} and the references therein). In this paper we will investigate the behaviour of \eqref{concaveconvex} as $q$ goes to $p$. The limit problem has the form

\begin{equation} \label{linearconvex}
\left\{
\begin{array}{r c l l}
-\Delta_p u  & = & \lambda\,u^{p-1}+ u^{r-1} &  \textrm{in
}\Omega, \\  u & = & 0 & \textrm{on }\partial\Omega,
\end{array}\right.
\end{equation}
and we will refer to it as \emph{linear-convex} problem. If we denote the first eigenvalue of the $p$-laplacian under Dirichlet boundary conditions by $\lambda_1(p)$, one can show that Problem \eqref{linearconvex} admits at least one positive solution for $\lambda < \lambda_1(p)$ (\cite[Section 3.3]{charroparini2}), and no positive solution for $\lambda \geq \lambda_1(p)$ (\cite[Proposition 3]{charroparini2}). Therefore, one can wonder whether the existence threshold $\Lambda_{q,r}$ tends to $\lambda_1(p)$ as $q \to p$. In this paper, we show that this is indeed the case, by proving lower and upper bounds for $\Lambda_{q,r}$ that are asymptotically optimal for $q \to p$.

The paper is structured as follows. After stating some preliminary results, in Section 3 we prove an upper bound for the existence threshold; this is done by obtaining a contradiction to an alternative definition of $\lambda_1(p)$. In Section 4 a lower bound for $\Lambda_{q,r}$ is proven by showing the existence of solutions to \eqref{concaveconvex} by means of an iteration method. The lower and the upper bound give the main result of this paper (Proposition \ref{mainresult}).

\section{Preliminary results}

Throughout the paper we will make use of the following classical regularity result.

\begin{prop}
Let $u \in W^{1,p}_0(\Omega)$ be a weak solution of
\begin{equation}\label{aabbcc}
\left\{
\begin{array}{r c l l}
-\Delta_p u  & = & f(u) &  \textrm{in
}\Omega, \\  u & = & 0 & \textrm{on }\partial\Omega,
\end{array}\right.
\end{equation}
where $f:\real\to \real$ is a Carath\'{e}odory function such that
\[ |f(u)|\leq C(1+|u|^q) \]
for some constant $C > 0$ and $q\in \left(1,\frac{np}{n-p}\right]$ if $p<n$, and $q \in (1,+\infty)$ otherwise. Then, $u \in L^\infty(\Omega)$ and therefore $u \in C^{1,\alpha}(\Omega)$ for some $\alpha \in (0,1)$.  
\end{prop}
\begin{proof}
The first part of the claim is a consequence of standard regularity estimates, while the second part follows from \cite{dibenedetto}.
\end{proof}

\subsection{The eigenvalue problem}

Let us consider the problem

\begin{equation} \label{eigenvalueproblem}
\left\{
\begin{array}{r c l l}
-\Delta_p u  & = & \lambda\,u^{p-1} &  \textrm{in
}\Omega, \\  u & = & 0 & \textrm{on }\partial\Omega.
\end{array}\right.
\end{equation}

We say that $\lambda \in \real \setminus \{0\}$ is an eigenvalue of the $p$-laplacian if there exists a nontrivial weak solution of \eqref{eigenvalueproblem}. The least eigenvalue can be found by minimizing the so-called Rayleigh quotient:

\[ \lambda_1(p):= \inf_{u \in W^{1,p}_0(\Omega)\setminus \{0\}} \frac{\int_\Omega |\nabla u|^p}{\int_\Omega |u|^p}. \]

It is known that the corresponding first eigenfunction is unique (up to a multiplicative constant), and does not change sign in $\Omega$; therefore it can be considered as strictly positive. Moreover, $\lambda_1(p)$ is isolated. The first eigenvalue can also be characterized as follows (see \cite{birindellidemengel}):

\begin{equation} \label{eigenvaluevaradhan}
\lambda_1(p) = \sup \{\lambda \geq 0\,|\,\exists v \in W^{1,p}_0(\Omega),\,v > 0 \text{ in }\Omega,\,-\Delta_p v \geq \lambda v^{p-1} \}. 
\end{equation}

\subsection{The linear-concave problem}

In the following we will state some preliminary results for the linear-concave problem

\begin{equation} \label{linearconcave}
\left\{
\begin{array}{r c l l}
-\Delta_p u  & = & \lambda\,u^{p-1}+ u^{q-1} &  \textrm{in
}\Omega, \\  u & = & 0 & \textrm{on }\partial\Omega,
\end{array}\right.
\end{equation}
where $1<q<p$. 

\begin{prop} \label{uniquenesslinearconcave} \cite[Lemma 4.1]{abdellaouiperal}  Let $u$, $v \in W^{1,p}_0(\Omega) \cap C^1(\Omega)$ be two strictly positive functions satisfying
\[ -\Delta_p u \leq \lambda u^{p-1} + u^{q-1}, \]
\[ -\Delta_p v \geq \lambda v^{p-1} + v^{q-1}. \] 
Then, $u \leq v$ in $\Omega$.
\end{prop}

\begin{prop} We have the following results:
\begin{enumerate}
\item[a)] For $\lambda \in [0,\lambda_1(p))$, there exists a unique positive solution of Problem \eqref{linearconcave}.
\item[b)] For $\lambda \in [\lambda_1(p), + \infty)$, Problem \eqref{linearconcave} admits no positive solution.
\end{enumerate}
\end{prop}
\begin{proof}
The proof of \emph{a)} can be found for instance in \cite[Section 4]{garciaperal}, and the proof of \emph{b)} in \cite[Proposition 3]{charroparini2}.
\end{proof}

\begin{prop} \label{upperbound}
For $\lambda \in [0,\lambda_1(p))$, let $u_{q,\lambda}$ be the positive solution of \eqref{linearconcave}. Then
\[ (\lambda_1(p)-\lambda)^{-\frac{1}{p-q}} \leq \|u_{q,\lambda}\|_\infty \leq (\|u_{q,0}\|_\infty^{q-p}-\lambda)^{-\frac{1}{p-q}}.\]
\end{prop}

\begin{proof} 
We know already that $u_{q,\lambda}\in L^\infty(\Omega)$; let $M=\|u_{q,\lambda}\|_\infty$ and define $v:=M^{-1}u_{q,\lambda}$. The function $v$ satisfies
\[ -\Delta_p v = \lambda v^{p-1} + M^{q-p} v^{q-1} \geq (\lambda + M^{q-p}) v^{p-1},\]
since $\|v\|_\infty = 1$. Hence, $v$ is a positive supersolution; by the characterization of $\lambda_1(p)$ in \eqref{eigenvaluevaradhan} it follows
\[ \lambda  + M^{q-p} \leq \lambda_1(p)\]
and hence
\[ (\lambda_1(p)-\lambda)^{-\frac{1}{p-q}} \leq \|u_{q,\lambda}\|_\infty .\]
Let us now prove the second part of the claim. It holds
\[ 1 = \|v\|_\infty \leq (\lambda + M^{q-p})^{\frac{1}{p-q}}\|u_{q,0}\|_\infty\]
and so
\[  \|u_{q,\lambda}\|_\infty \leq (\|u_{q,0}\|_\infty^{q-p}-\lambda)^{-\frac{1}{p-q}}.\]
\end{proof}

\begin{prop}
For $\lambda \in [0,\lambda_1(p))$, let $u_{q,\lambda}$ be a solution of \eqref{linearconcave}. Then
\[ \|u_{q,\lambda}\|_\infty^{q-p} \to \lambda_1(p)-\lambda \]
as $q \to p$.
\end{prop}
\begin{proof}
By \cite[Corollary 3.3]{ercole} we know that
\[ \|u_{q,0}\|_\infty^{q-p} \to \lambda_1(p) \]
as $q \to p$. The claim then follows from Proposition \ref{upperbound}.
\end{proof}

\noindent This implies in particular that the quantity
\begin{equation} \label{definitionofc}
c(q,\lambda):= \frac{\|u_{q,\lambda}\|_\infty}{(\lambda_1(p)-\lambda)^{-\frac{1}{p-q}}}
\end{equation}
is such that $c(q,\lambda)^{p-q} \to 1$ as $q \to p$. 

\subsection{The iteration method} \label{iterationmethod}

We will detail here for later reference the construction of a viscosity solution to Problem \eqref{concaveconvex} by iteration. We point out that the notions of continuous weak solutions and viscosity solutions are equivalent for problems \eqref{concaveconvex} and \eqref{linearconvex} in the whole range $p>1$. The proof of this fact is a modification of the arguments in \cite{Julin-Juutinen}.

First, we will assume the existence of a subsolution $\underline{u}$ and a supersolution  $\overline{u}$ of  \eqref{concaveconvex}
such that $\underline{u}=\overline{u}=0$ on $\partial\Omega$ and 
 $0<\underline{u}\leq\overline{u}$ in $\Omega$.

Let  $w_1(x)$ be the viscosity solution of
\[
\left\{
\begin{split}
-&\Delta_p w_1=\Lambda\,\underline{u}^{q-1}+\underline{u}^{r-1}\quad\text{in}\ \Omega,\\
&w_1=0\quad\text{on}\ \partial\Omega.
\end{split}
\right.
\]
Such a solution $w_1$ exists since $\underline{u}$ is a subsolution and $\overline{u}$ is a supersolution, as
\[
-\Delta_p\overline{u}\geq\Lambda\overline{u}^{q-1}+\overline{u}^{r-1}\geq\Lambda\underline{u}^{q-1}+\underline{u}^{r-1}.
\]
By comparison and the Perron method we get that there exists a unique $w_1$ such that
\[
\underline{u}\leq w_1\leq\overline{u}\quad \text{in}\ \Omega.
\]
Then, define $w_2$, solution of
\[
\left\{
\begin{split}
-&\Delta_pw_2=\Lambda\,w_1^{q-1}+w_1^{r-1}\quad\text{in}\ \Omega,\\
&w_2=0\quad\text{on}\ \partial\Omega.
\end{split}
\right.
\]
In this case, $w_1$ is a subsolution and $\overline{u}$ a supersolution, since by monotonicity
\[
-\Delta_p w_1=\Lambda\,\underline{u}^{q-1}+\underline{u}^{r-1}\leq \Lambda\,w_1^{q-1}+w_1^{r-1},
\]
while
\[
-\Delta_p\overline{u}\geq\Lambda\,\overline{u}^{q-1}+\overline{u}^{r-1}\geq\Lambda\,w_1^{q-1}+w_1^{r-1}.
\]
We have $w_1=\overline{u}=0$ on $\partial\Omega$, and hence by comparison
and the Perron method we get the existence of $w_2$ such that
\[
\underline{u}\leq w_1\leq w_2\leq\overline{u}\quad\text{in}\ \Omega.
\]

Iterating this procedure we can construct a sequence  $\{w_k\}_{k\geq1}$ of solutions of
\begin{equation}\label{eq:exSI:eq:22}
\left\{
\begin{split}
-&\Delta_pw_k=\Lambda\,w_{k-1}^{q-1}+w_{k-1}^{r-1}\quad\text{in}\ \Omega,\\
&w_k=0\quad\text{on}\ \partial\Omega.
\end{split}
\right.
\end{equation}
such that
\[
\underline{u}\leq w_1\leq w_2\leq\ldots \leq w_{k-1}\leq w_k\leq\overline{u}\quad\text{in}\ \Omega.
\]
In particular, for every $x\in\Omega$, the sequence $\{w_k(x)\}_{k\geq1}$ is bounded and is non-decreasing, hence convergent. We denote  $u(x)$ the pointwise limit of the  $w_k$. Then, there exists a subsequence $k'\to\infty$ such that
\[
w_{k'} \to u\quad\text{uniformly in }\overline{\Omega}.
\]

Since the sequence $w_k$ is monotonically increasing, the whole sequence converges uniformly to  $u$. It is then easy to prove that $u$ is a viscosity solution of \eqref{concaveconvex}.

\section{Upper bound for the existence threshold}

In this section we give an explicit value $\hat{\Lambda}_p>0$ such that no positive weak solution of \eqref{concaveconvex} exists for $\Lambda>\hat{\Lambda}_p$. Therefore, it is clear that $\Lambda_{q,r} \leq \hat{\Lambda}_p$.

\begin{prop} \label{upperboundforthreshold}
Let 
\[ \hat{\Lambda}(p,q,r) = \lambda_1(p)^{\frac{r-q}{r-p}}(r-p)\left(\frac{(p-q)^{p-q}}{(r-q)^{r-q}} \right)^{\frac{1}{r-p}}.\]
Then Problem \eqref{concaveconvex} does not admit any positive solution for $\Lambda >  \hat{\Lambda}_p$. As a consequence, $\Lambda_{q,r} \leq \hat{\Lambda}(p,q,r)$ and therefore
\[ \limsup_{q \to p} \Lambda_{q,r} \leq \lambda_1(p). \]
\end{prop}
\begin{proof}
Suppose by contradiction that Problem  \eqref{concaveconvex} has a positive solution  $u_{\Lambda}\in W_0^{1,p}(\Omega)$ for some $\Lambda>\hat\Lambda(p,q,r)$. Then, we can choose $\epsilon>0$ small enough such that 
\[
\Lambda>\mu^\frac{r-q}{r-p}(r-p)\left(\frac{(p-q)^{p-q}}{(r-q)^{r-q}}\right)^\frac{1}{r-p}.
\]
for $\mu=\lambda_1(p)+\epsilon$.

We claim that then
\begin{equation}\label{eq:CCNO:supersolpp}
-\Delta_{p}u_{\Lambda}=\Lambda\,u_{\Lambda}^{q-1}+u_{\Lambda}^{r-1}>\mu u_{\Lambda}^{p-1} \quad\text{in}\ \Omega,
\end{equation}
in the weak sense, a contradiction to Proposition \ref{eigenvaluevaradhan}. 
In order to prove the claim, it is enough to see that
\[
\min_{t>0}\Phi_\Lambda(t)>\mu\quad\text{where}\quad\Phi_\Lambda(t)=\Lambda\,t^{q-p}+t^{r-p}.
\]
It is elementary to check that
\[
\frac{d}{dt}\Phi_\Lambda(t)=0\quad\Leftrightarrow\quad t_\Lambda=\left(\frac{\Lambda\big(p-q\big)}{\big(r-p\big)}\right)^\frac{1}{r-q},
\]
which is a minimum. As  $\Phi_\Lambda(t)\to\infty$ when $t\to0$ and  $t\to\infty$, it is a global minimum. Then,
\[
\min_{t>0}\Phi_\Lambda(t)=\Phi_\Lambda(t_\Lambda)= \frac{\Lambda^\frac{r-p}{r-q}(r-q)}{(p-q)^\frac{p-q}{r-q}(r-p)^\frac{r-p}{r-q}}>\mu
\]
because of our choice of  $\Lambda$. The first part of the proposition is therefore proved, while the second part is an easy consequence of the definition of $\Lambda_{q,r}$.

\end{proof}

\section{Lower bound for the existence threshold}\label{section.lower.bound}

The aim is to prove that Problem \eqref{concaveconvex} admits a solution for $\Lambda\leq \widetilde{\Lambda}(p,q,r)$, where the value of $\widetilde{\Lambda}(p,q,r)$ is given explicitly as a function of $p$, $q$, $r$, $\lambda_1(p)$ and the quantity $c(q,\lambda)$ defined in \eqref{definitionofc}. This clearly implies that $\widetilde{\Lambda}(p,q,r) \leq \Lambda_{q,r}$.
\begin{prop}\label{prop.existence.small.lambda}
Let
\begin{equation}
\label{tildelambda}
 \widetilde{\Lambda}(p,q,r) = \lambda_1(p)^{\frac{r-q}{r-p}} (r-p) \left(\frac{(p-q)^{p-q}}{(r-q)^{r-q}} \right)^{\frac{1}{r-p}} c(q,t_q\lambda_1(p))^{q-p}, 
\end{equation}
where $t_q=\frac{p-q}{r-q}$, and $c(q,\lambda)$ is the constant defined in \eqref{definitionofc}. Then Problem \eqref{concaveconvex} admits a positive solution for every $\Lambda \in [0, \widetilde{\Lambda}(p,q,r)]$.
\end{prop}
\begin{proof}
In order to find a solution of Problem \eqref{concaveconvex}, we will use the iteration method described in Section \ref{iterationmethod}. Therefore, we need to find a subsolution $\underline{u}$ and a supersolution $\overline{u}$ such that $\underline{u}\leq \overline{u}$. A subsolution of Problem \eqref{concaveconvex} is given by the positive solution of
\[
\left\{
\begin{array}{r c l l}
-\Delta_p w  & = &  \Lambda w^{q-1} &  \textrm{in
}\Omega \\  w & = & 0 & \textrm{on }\partial\Omega.
\end{array}\right.
\]
Let us now look for a supersolution. Let $u_{\Lambda}$ be the positive solution of
\[
\left\{
\begin{array}{r c l l}
-\Delta_p u  & = & \lambda\,u^{p-1}+ \Lambda u^{q-1} &  \textrm{in
}\Omega \\  u & = & 0 & \textrm{on }\partial\Omega,
\end{array}\right.
\]
where $\lambda \in [0,\lambda_1(p))$ will be chosen later. By scaling, one notices that $u_{\Lambda}=\Lambda^{\frac{1}{p-q}} u_1$, and therefore
\[ \|u_{\Lambda}\|_\infty^{p-q} =  \frac{\Lambda \cdot c(q,\lambda)^{p-q}}{\lambda_1(p)-\lambda} \]
by \eqref{definitionofc}.
In order for $u_{\Lambda}$ to be a supersolution, it is necessary and sufficient that
\[ \lambda\,u_{\Lambda}^{p-1}+ \Lambda u_{\Lambda}^{q-1} \geq \Lambda\,u_{\Lambda}^{q-1}+ u_{\Lambda}^{r-1}.\]
This is equivalent to
\[ \lambda \geq u_{\Lambda}^{r-p}(x) \]
for every $x\in\Omega$, which is equivalent to
\[ \lambda \geq \|u_{\Lambda}\|_\infty^{r-p}. \]
This is in turn equivalent to
\[ \lambda^\frac{p-q}{r-p} \geq  \frac{\Lambda \cdot c(q,\lambda)^{p-q}}{\lambda_1(p)-\lambda}\]
which is equivalent to
\[ \Lambda \leq  \lambda^\frac{p-q}{r-p} (\lambda_1(p)-\lambda) \cdot c(q,\lambda)^{q-p}.\]
Set now $\lambda = t_q \lambda_1(p)$ with $t_q=\frac{p-q}{r-q}\in [0,1)$, which implies $\Lambda\leq\tilde\Lambda$ (the value $t_q$ is chosen since it is the maximal value for the function
$f(t)= t^\frac{p-q}{r-p} (1-t)$).
Notice that, if we choose $\lambda = t_q \lambda_1(p)$, we obtain a supersolution every $\Lambda\leq\tilde\Lambda$. Since $w \leq u_{\Lambda}$ by Proposition \ref{uniquenesslinearconcave}, it is possible to apply the iteration method in order to find a solution of \eqref{concaveconvex}.
 
\end{proof}

\begin{prop} \label{lowerboundforthreshold}
We have that
\[ \Lambda_{q,r} \geq \lambda_1(p)^{\frac{r-q}{r-p}}(r-p)\left(\frac{(p-q)^{p-q}}{(r-q)^{r-q}} \right)^{\frac{1}{r-p}} \cdot c(q,t_q\lambda_1(p))^{q-p} \]
and therefore
\[ \liminf_{q \to p} \Lambda_{q,r} \geq \lambda_1(p). \] 
\end{prop}
\begin{proof}
The first inequality is an easy consequence of Proposition \ref{prop.existence.small.lambda}. Reasoning as in Proposition \ref{upperbound}, we have that
\[ (\lambda_1(p)-t_q\lambda_1(p)) \geq \|u_{q,t_q\lambda_1(p)}\|^{q-p}_\infty \geq (\|u_{q,0}\|_\infty^{q-p}-t_q\lambda_1(p))\]
and therefore
\[ c(q,t_q\lambda_1(p))^{q-p} \to 1 \]
as $q\to p$. This implies the second part of the claim.
\end{proof}

We can finally state the main result of this paper.

\begin{prop} \label{mainresult} Let $\Lambda_{q,r}$ be the existence threshold for Problem \eqref{concaveconvex}. Then,
 \[ \lim_{q \to p} \Lambda_{q,r} = \lambda_1(p). \] 
\end{prop}
\begin{proof}
The proof is a consequence of Propositions \ref{upperboundforthreshold} and \ref{lowerboundforthreshold}
\end{proof}

\bibliographystyle{amsplain}
\bibliography{bibliofernandoeneathird}

\end {document}